\documentclass[twocolumn,prd,aps,amsfonts,amsmath,mathtools,amssymb,showpacs,preprintnumbers,nofootinbib,superscriptaddress,noshowpacs]{revtex4}
\usepackage{tikz}
\usetikzlibrary{patterns}
\usepackage{mathtext}
\usepackage[T1]{fontenc}
\usepackage[english,russian]{babel}
\usepackage{tabularx}
\usepackage{longtable}
\newcolumntype{C}[1]{>{\centering\arraybackslash}b{#1}}
\newcolumntype{D}[1]{>{\centering\arraybackslash}b{#1}}
\newcolumntype{R}[1]{>{\raggedleft\arraybackslash}b{#1}}

\def\e{\mathrm{e}}

\def\I{\mathrm{i}}

\def\Z{\mathbb{Z}}

\newcommand{\be}{\begin{eqnarray}}
\newcommand{\ee}{\end{eqnarray}}
\newcommand{\nn}{\nonumber}

\newcommand{\reals}{\mathbb{R}}
\newcommand{\ints}{\mathbb{Z}}
\newcommand{\oct}{\mathbb{O}}

\renewcommand{\d}{{\mathrm d}}

\newcommand{\R}{\mathbb{R}}

\newcommand{\mak}[1]{\mathfrak{#1}}
\def\Lo{\CYRL}
\DeclareMathOperator{\Li}{Li}
\renewcommand{\C}{\mathbb{C}}

\renewcommand\Re{\operatorname{Re}}
\renewcommand\Im{\operatorname{Im}}
\def\vol{\text{vol}}
\newcommand{\braket}[2]{\langle #1 | #2 \rangle}

\newcommand\cF{{\cal{F}}}
\newcommand\cH{{\cal{H}}}
\newcommand\cC{{\cal{C}}}

\newcommand{\bu}{{\mathbf{u}}}
\newcommand{\bX}{{\mathbf{x}}}
\newcommand{\ba}{{\mathbf{a}}}
\newcommand{\beb}{{\mathbf{e}}}
\newcommand{\bl}{{\boldsymbol{\lambda}}}
\newcommand{\bt}{{\boldsymbol{\theta}}}

\newcommand\de{\delta}
\newcommand\bd{\bar\delta}

\newcommand\mhlinee{\hline & & & & & & & \\[-2mm]}
\newcommand\mhlines{\hline & & & & & & \\[-2mm]}

\begin{document}
\selectlanguage{english}

\preprint{AEI-2011-012}

\title{On fundamental domains and volumes of hyperbolic Coxeter-Weyl groups}
\author{Philipp Fleig}
\affiliation{Max-Planck-Institut f\"ur Gravitationsphysik, Albert-Einstein-Institut, 
Am M\"uhlenberg 1, 14476 Potsdam, Germany}
\affiliation{Universit\'e de Nice-Sophia Antipolis, Parc Valrose, 06108 Nice Cedex 2,
France}
\author{Michael Koehn}
\author{Hermann Nicolai}
\affiliation{Max-Planck-Institut f\"ur Gravitationsphysik, Albert-Einstein-Institut, 
Am M\"uhlenberg 1, 14476 Potsdam, Germany}

\date{November 9, 2011}

\begin{abstract}
We present a simple method for determining the shape of fundamental domains of 
generalized modular groups related to Weyl groups of hyperbolic Kac-Moody 
algebras. These domains are given as subsets of certain generalized upper half planes,
on which the Weyl groups act via generalized modular transformations.
Our construction only requires the Cartan matrix of the underlying finite-dimensional 
Lie algebra and the associated Coxeter labels as input information. We present 
a simple formula for determining the volume of these fundamental 
domains. This allows us to re-produce in a simple manner the known values for 
these volumes previously obtained by other methods.
\end{abstract}
\maketitle

\section{Introduction}

Constructions of fundamental domains of generalized modular groups usually 
rely on geometric considerations. By considering the different 
possible symmetry transformations acting on some generalized upper-half plane, 
the precise shape of the fundamental domain is narrowed down step-by-step until 
one arrives at its final shape. Especially for higher rank groups (such
as $SL_n(\ints)$) this poses a considerable computational and combinatorial 
problem since one has to consider a large number of possible successive symmetry transformations (already the determination of the fundamental domain of the standard 
modular group $PSL_2(\ints)$ along these lines takes more than two pages of computations, 
see e.g. \cite{Apo90}). Although one can show that the precise shape of the fundamental 
domain can be determined within a finite number of steps, in the actual computation of a 
domain it is not always clear how many steps are actually necessary.

In this paper we show that, at least for modular groups arising as (even) Weyl groups
of certain hyperbolic Kac--Moody algebras, such cumbersome constructions can be
altogether avoided. More specifically, we present an easy method for obtaining the 
complete geometric information about the associated fundamental domains. All we 
require as information for determining the explicit shape and volume is the Cartan 
matrix of the corresponding Kac-Moody algebra and its Coxeter labels. As we will  
demonstrate this construction works for {\em all} hyperbolic  Kac--Moody algebras 
\footnote{An indefinite Kac--Moody algebra is called {\em hyperbolic} if the removal of
 any one node from its Dynkin diagram leaves an algebra which is either affine
 or finite \cite{Kac90}.} 
$\mak{g}^{++}$ of over-extended type, which are generally obtained by extending a 
given finite dimensional simple Lie algebra $\mak{g}$ via its affine extension 
$\mak{g}^+$ by adding two nodes to the Dynkin diagram in a specified way. Likewise
it applies to the {\em twisted} algebras obtained by inverting the arrows in the Dynkin
diagram, because their Weyl groups are the same (but note that these twisted algebras, 
while being indefinite Kac--Moody algebras, in general are not of over-extended type). 
In particular, our construction also applies to those hyperbolic Kac--Moody algebras
whose even Weyl groups can be identified with generalized modular groups 
defined over rings of integers in division algebras \cite{FKN09}. The first 
example of such an identification was given in \cite{FF} where it was shown that the 
rank-3 hyperbolic Kac--Moody algebra $A_1^{++}$ (also denoted $AE_3$ or $\cF$ in 
the literature) has the usual modular group $PSL_2(\ints)$ as its even Weyl group, 
the full Weyl group being $W(A_1^{++}) = PGL_2(\ints)$. In \cite{FKN09} more
complicated examples were given, involving for instance the quaternionic integers
(Hurwitz numbers), and admitting a M\"obius-like realization \cite{KNP10}.
The most interesting (and most complicated) example is the even Weyl group 
$W^+(E_{10})$  which can be identified with the arithmetic group $PSL_2(\mathtt{O})$ 
(where $\mathtt{O}$ are octonionic integers, also called {\it octavians}). For this 
example we will explicitly display the coordinates of the vertices of the 
fundamental domain of the Weyl group.

Knowledge of the shape of the fundamental domain allows one to compute its volume.
In the non-linear realization of the hyperbolic Weyl group on some generalized upper 
half plane \cite{KNP10} (a hyperbolic space of constant negative curvature) the fundamental domains are realized as  higher dimensional simplices. We present a very simple general 
formula for the volume of the domain in terms of integrals involving a quadratic form 
which contains all the information about the Lie algebra $\mak g^{++}$ (see (\ref{main}) 
below).  We note that our considerations would also apply to cases where analogs of 
the so-called congruence subgroups of $PSL_2(\ints)$ can de defined: the volume 
is then simply a multiple of the original volume, with the factor equal to the index of the congruence subgroup in the given generalized modular group. Such congruence subgroups
presumably do exist for the generalized arithmetic groups studied in \cite{FKN09}, 
but we are not aware of any concrete results along these lines. 

As an historic aside, we mention that
the first computation of hyperbolic volumes in terms of the dihedral angles of 
the simplex under consideration is due to one of the inventors of hyperbolic
geometry, N.I. Lobachevsky \cite{Lob04}. His results were extended by Schl\"afli 
and Coxeter \cite{Cox35}, see also Vinberg \cite{Vin93b}. Further work on this problem 
can be found in \cite{JKRT99} which gives a list of numerical values for the volumes of 
hyperbolic Coxeter simplices, as well as analytical expressions for some special cases. 
Using (\ref{main}) these values can be easily reproduced. We also note that in the physical context, these Coxeter simplices appear in the cosmological billiards setting, see \cite{KKN09,Koe11} for the implications of the quantum treatment of the cosmological billiards for an initial spacelike singularity.

\section{Hyperbolic roots and weights}

Let $\mak{g}$ be a finite-dimensional Lie algebra. We denote the simple roots 
of $\mak g$ by $\ba_i\in\R^n$ and their associated fundamental weights by $\bl_i$, 
where $i=1,\ldots ,n$ with $n=\mathrm{Rank}(\mak g)$ (see e.g. \cite{Hum72} for
details). With the Cartan matrix of $\mak g$ 
\be
A_{ij} = \langle \ba_i | \ba_j\rangle  \equiv \frac{2 \ba_i \cdot \ba_j}{\ba_j\cdot \ba_j} 
\ee
we define the {\em symmetrized Cartan matrix} $B_{ij}$ as
\be
B_{ij} \equiv (AD)_{ij} =2 \ba_i \cdot \ba_j = A_{ij} \, \ba_j^2 \ , \label{symmcart}
\ee
where $\ba_j^2\equiv \ba_j\cdot \ba_j$ and there is no summation 
over double indices. Unlike $A_{ij}$, the matrix $B_{ij}$ and the symmetrizing 
matrix $D_{ij} = \delta_{ij} \ba_j^2$ depend on the normalization of $\ba_j$. Following \cite{FKN09} we choose this normalization such that always $\bt^2=1$ for
 the highest root
\be
\bt = \sum_{j=1}^n m_j \ba_j
\ee
with the Coxeter labels $m_j$. When $\bt$ is a {\em long} root
we therefore have $\ba_j^2=1$ for the long roots. 

The associated fundamental weights $\bl_j$ constitute a basis dual to the
simple roots \cite{Hum72}
\be\label{lambda}
\braket{\bl_i}{\ba_j}\equiv \frac{2\bl_i \cdot \ba_j}{\ba_j\cdot \ba_j}=\delta_{ij}
\ee
implying
\be
\bl_i\cdot \ba_j = \frac{1}{2} \, \delta_{ij} \, \ba_j^2 \ .
\ee
With the inverse Cartan matrix $A^{-1}$ we thus have
\be
\bl_i = \sum_k (A^{-1})_{ik} \ba_k \ ,
\ee
from which we deduce
\be
\bl_i \cdot \bl_j = \frac{1}{2} (A^{-1})_{ij} \, \ba_j^2 \ .
\ee
or
\be\label{lambda1}
\bl_i \cdot \bl_j = \frac{1}{2} \ba_i^2 \,(B^{-1})_{ij} \, \ba_j^2 \ .
\ee
Next we consider the {\em hyperbolic extension} $\mak g^{++}$ of the finite-dimensional 
algebra $\mak g$ obtained by adjoining to the Dynkin diagram of $\mak g$  the affine node (labeled `0') and the over-extended node (labeled `$-1$'). This entails extending 
the Euclidean root space $\R^n$ to the {\em Lorentzian} space $\R^{1,n+1}=\R^{1,1} \oplus \R^n$. 
We denote the roots of $\mak g^{++}$ by $\alpha_I$, $I=-1,0,1,\ldots , n$, and define 
them according to
\be
\alpha_{-1}\equiv -\de - \bd \ , \quad \alpha_0\equiv \de- \bt \ , \quad \alpha_i\equiv \ba_i
\ee
with the affine null vector $\de\in\R^{1,1}$ and the conjugate null vector $\bd\in\R^{1,1}$ 
obeying $\de \cdot \bd=\frac{1}{2}$. In this way we obtain the Cartan matrix of 
${\mak g}^{++}$ as
\be
A_{IJ}= \langle \alpha_I | \alpha_J\rangle \equiv 
     \frac{2 \alpha_I \cdot \alpha_J}{\alpha_J \cdot \alpha_J}
\ee
with the {\em Lorentzian} inner product 
\be
\alpha_I\cdot \alpha_J \equiv \eta_{\mu\nu} \alpha_I^\mu \alpha_J^\nu 
 \ee
where the signature of $\eta_{\mu\nu}$ is  $(-+\cdots +)$.
Notice that the affine and over-extended simple roots are also
normalized as $\alpha_{-1}^2 = \alpha_0^2 =1$. The normalization
$\bt^2=1$ is necessary to obtain a single line between the affine
and the hyperbolic node (connecting $\alpha_0$ and $\alpha_{-1}$).

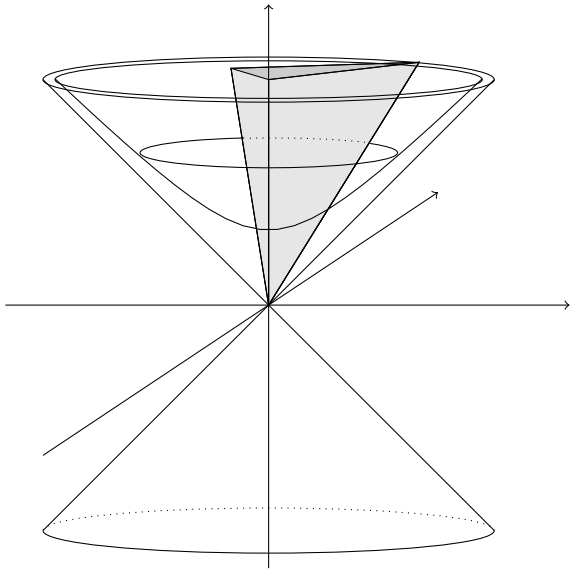
\begin{figure}
\begin{center}
\begin{tikzpicture}[domain=-2.84:2.84]
\filldraw[fill=gray!40!white, draw=black] (0,0) -- (-0.5,3.15) -- (2,3.23) -- (0,0);
\filldraw[fill=gray!20!white, draw=black] (0,0) -- (0.0,3.0) -- (-0.5,3.15) -- (0,0);
\filldraw[fill=gray!20!white, draw=black] (0,0) -- (0.0,3.0) -- (2,3.23) -- (0,0);
\draw (0,0) -- (-0.5,3.15) -- (2,3.23) -- (0,0);
\draw (0,0) -- (0.0,3.0) -- (2,3.23) -- (0,0);
\draw (-3.0,-3.0) -- (3.0,3.0);
\draw (-3.0,3.0) -- (3.0,-3.0);
\draw (0.0,3.0) ellipse (3 and 0.3);
\draw[color=black] (0.0,3.0) ellipse (2.84 and 0.25);
\draw[color=black] (1.35,2.15) arc (38:-259:1.71 and 0.2);
\draw[dotted,color=black] (-0.25,2.22) arc (101:43:1.71 and 0.2);
\draw[dotted] (3.0,-3.0) arc (0:180:3 and 0.3);
\draw (-3.0,-3.0) arc (180:360:3 and 0.3);
\draw[color=black]	plot ( \x ,{sqrt(\x*\x+1.0)});
\draw[->] (0,-3.5) -- (0,4);
\draw[->] (-3.5,0) -- (4,0);
\draw[->] (-3,-2) -- (2.25,1.5);
\end{tikzpicture}
\caption{Sketch of the fundamental Weyl chamber $\cC_0$ as a wedge inside the forward light cone that is 
 intersected by the unit hyperboloid.}\label{fig:kegel}
\end{center}
\end{figure}

The fundamental weights $\Lambda_I$ for the hyperbolic extension
$\mak g^{++}$ are defined in analogy with (\ref{lambda}) 
\be\label{scalarhyp}
\braket{\Lambda_I}{\alpha_J}\equiv \frac{2\Lambda_I\cdot \alpha_J}{\alpha_J \cdot \alpha_J} = \delta_{IJ}\ .
\ee
By a standard construction (see e.g. \cite{DHJN01}), the fundamental weights of
$\mak g^{++}$ can be expressed in terms of the null vectors $\de$ and $\bd$ 
and the finite weights $\bl_j$ as
\be\label{Lambda}
\Lambda_{-1}=-\de \ , \quad \Lambda_0=\bd - \de \ , \quad \Lambda_j=n_j \Lambda_0+\bl_i \ .
\ee
The coefficients $n_j$ are fixed by requiring $\alpha_0\cdot \Lambda_j=0$ 
(cf. \eqref{scalarhyp}), which gives
\be
n_j =  m_j \ba_j^2 \ ,
\ee
The fundamental Weyl chamber $\cC_0 \subset \R^{1,n+1}$ is 
\be\label{F1}
\cC_0 := \big\{ X\in\R^{1,n+1} \, | \, X \cdot \alpha_I \geq 0 \;\mbox{ for $I=-1,0,1,...,n$}  \big\}
\nonumber
\ee
With the fundamental weights $\Lambda_I$ one obtains a more convenient
representation of $\cC_0$
\be\label{F2}
\cC_0 = \big\{ X \in \R^{1,n+1}\, |\,  X = \sum_I s_I \Lambda_I \;\; \mbox{with $s_I\geq 0$ for
                     all $I$} \big\}\nn\\
\ee
The null vector $\de$ lies on the forward 
light-cone in root space. The fundamental Weyl chamber $\cC_0$ is the convex hull of the hyperplanes 
orthogonal to the simple roots of the algebra. The fundamental weights are vectors 
pointing along the edges of $\cC_0$. In other words, $\cC_0$ is a `wedge'  in
$\R^{1,n+1}$. For the  hyperbolic algebras ${\mak g}^{++}$ of over-extended type
considered here this wedge lies inside the forward lightcone,  always touching it with the
lightlike weight vector $\Lambda_{-1}$, while all other fundamental weights obey
$\Lambda^2_j\leq 0$. By contrast, for general indefinite (Lorentzian) ${\mak g}^{++}$ 
the fundamental Weyl chamber may stretch beyond the lightcone and also contain space-like vectors.
A schematic picture of the fundamental Weyl chamber $\cC_0$ for hyperbolic $\mak{g}^{++}$ is  
shown in Fig.\,\ref{fig:kegel}. We have included the forward light-cone and the intersecting unit hyperboloid.

As it turns out the assumptions made suffice to cover all cases of interest.
This concerns in particular the {\em twisted} algebras: as these are obtained
by inverting the arrows in the relevant Dynkin diagrams, the associated
Coxeter Weyl groups, not being sensitive to the direction of the arrows,
coincide with those of the untwisted diagrams. We therefore note the following 
isomorphisms of Weyl groups using Kac' notation \cite{Kac90}: 
\be
W(G_2^{(1)+}) &\cong& W(D_4^{(3)+} )  \nn\\
W(B_n^{(1)+}) &\cong& W(A_{2n-1}^{(2)+} )  \nn\\
W(C_n^{(1)+}) &\cong& W(D_{n+1}^{(2)+})  \nn\\
W(F_4^{(1)+}) &\cong& W(E_{6}^{(2)+})   
\ee
where the superscript $^+$ on the r.h.s. indicates the extension
of the affine algebra by another node. But note that the twisted
algebras, though perfectly well-defined as indefinite Kac--Moody
algebras, are not necessarily of over-extended type.
In the notation of Fuchs
and Schweigert \cite{FS}, the later three isomorphisms are
\be
W(B_n^{(1)+}) &\cong&  W(C_{n}^{(2)+}) \nn\\
W(C_n^{(1)+}) &\cong& W(B_{n}^{(2)+}) \nn\\
W(F_4^{(1)+}) &\cong&   W(F_{4}^{(2)+}) 
\ee
The corresponding volumes of the fundamental domains
therefore also coincide.

\section{Volume Formula}

The linear action of the Weyl group in $\R^{1,n+1}$ preserves the (Lorentzian) length,
and therefore induces a non-linear {\em modular action} on the forward unit hyperboloid 
\be
X\cdot  X \equiv -x^+x^-+ \bX^2=  - 1 \;\; , \quad x^{\pm}>0
\ee
with light-cone coordinates $x^\pm\equiv (x^0\pm x^{n+1})/\sqrt{2}$ in $\R^{1,1}$ 
and $\bX\in\R^n$. For the cases $n=1,2,4$ and 8 studied in \cite{FKN09}, where the
dual of the Cartan subalgebra of $\mak g$ can be endowed with the structure of a division
algebra, the induced non-linear action takes the form of a generalized M\"obius
transformation over a (possibly non-commutative and non-associative)
ring of integers.

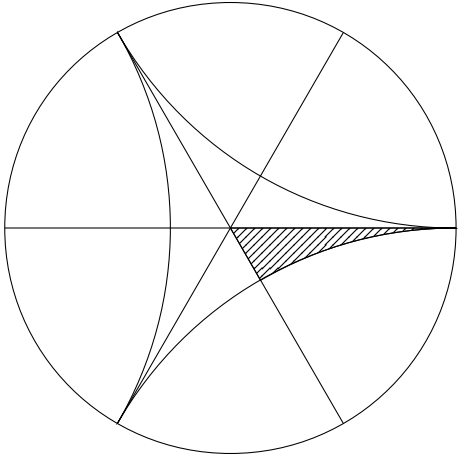
\begin{figure}
\begin{center}
\begin{tikzpicture}
\draw (0,0) circle (3);
\filldraw[fill=gray!20!white, draw=black,pattern=north east lines] (-3,0) -- (3,0) arc (90:120:5.2) -- (0,0);
\draw (-1.5,-2.598) -- (1.5,2.598);
\draw (-1.5,2.598) -- (1.5,-2.598);
\draw (3,0) arc (90:150:5.2) -- (-1.5,-2.598) arc (-30:30:5.2) -- (-1.5,2.598) arc (-150:-90:5.2) -- (3,0);
\draw (0,0) -- (3,0);
\end{tikzpicture}
\caption{Example of a fundamental domain on the Poincar\'e disk obtained by intersecting
the Weyl chamber with the (compactified) unit hyperboloid, here for the 
algebra $A_1^{++}$.}\label{fig:poincare}
\end{center}
\end{figure}

\begin{figure}
\begin{center}
\begin{tikzpicture}
\draw[dashed] (0,0) -- (0,1.2);
\draw[->] (0,2.8) -- (0,6) node[right]{$v$};
\draw (3,0) arc (0:180:3);
\draw (-3,0) arc (-180:0:3 and 1);
\draw[dotted] (-3,0) arc (180:0:3 and 1);
\filldraw[fill=gray!20!white, draw=black, pattern=north east lines] (0,2.8) arc (110:166:2.05 and 3) -- (-1.3,0.7) arc (121.5:58.5:2.5 and 3) -- (1.3,0.7) arc (14:70:2.05 and 3);
\draw (-1.3,0.7) -- (-1.3,5.5);
\draw[dashed] (1.3,5.5) -- (1.3,6.5);
\draw[dashed] (-1.3,5.5) -- (-1.3,6.5);
\draw[dashed] (0,6) -- (0,7);
\draw (1.3,0.7) -- (1.3,5.5);\draw[dashed] (-1.3,0.7) -- (-1.3,-0.5);
\draw[dashed] (1.3,0.7) -- (1.3,-0.5);
\filldraw[fill=gray!20!white, draw=black, pattern=north east lines] (-1.3,-0.5) -- (0,0) -- (1.3,-0.5) -- (-1.3,-0.5);
\end{tikzpicture}
\caption{Schematic depiction of a Weyl chamber on the UHP, corresponding 
in this case to $A_2^{++}$.}\label{fig:lqs}
\end{center}
\end{figure}
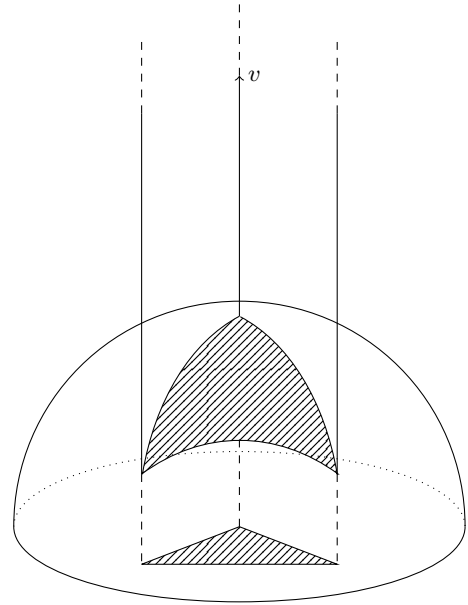

The intersection of the fundamental Weyl chamber $\cC_0$ with the 
unit hyperboloid defines a corresponding fundamental domain $\cF_0$ 
on the unit hyperboloid.  The corresponding domain on the (compactified)
unit hyperboloid ({\it alias} the Poincar\'e disk) is depicted in 
Fig.\,\ref{fig:poincare}. In the remainder, however,  we will study this domain 
as a subset of the {\em generalized (Poincar\'e) upper half plane} $\cH$
rather than the unit hyperboloid \footnote{Note that the fundamental domain
   $\cF_0$ is half of the fundamental domain $\cF$ of the ordinary
   modular group. The latter corresponds to the {\em even} subgroup
   of the Weyl group.}.  This upper half plane is defined as
\be
\cH \equiv \cH_{n+1} := \big\{ (\bu, v)\, | \, \bu\in \R^n \,,\, v >0 \big\}
\ee
and is thus of dimension $n+1$. $\cH_{n+1}$ is isometric to the forward unit hyperboloid
in $\reals^{1,n+1}$ by means of the standard coordinate transformation
\be\label{UHP1}
x^- = \frac1{v} \;, \quad    x^+ = v + \frac{\bu^2}v \;, \quad \bX= \frac{\bu}{v} 
\ee
The Minkowskian line element is transformed to 
\be
\d s^2 = \frac{\d\bu^2 + \d v^2}{v^2}
\ee
where, of course, $\d\bu^2 \equiv \d u_1^2 + \cdots + \d u_n^2$.
The fundamental domain $\cF_0\subset\cH$ is now rather easy to determine
from the representation (\ref{F2}) by identifying the points where the rays
along $\Lambda_I$ `pierce' the unit hyperboloid, and then mapping these points
to $\cH$ by means of (\ref{UHP1}). We first notice that the over-extended fundamental 
weight $\Lambda_{-1}$ ({\it alias} the affine null vector $\de$) corresponds 
to the `cusp' at infinity in $\cH$ with coordinates $v=\infty\,,\, \bu=0$ \cite{KNP10}, 
while $\Lambda_0$ corresponds to the point $v=1\,,\, \bu=0$ in $\cH$. From (\ref{Lambda}) we see that the remaining fundamental weights are mapped to the points
\be\label{UHPMap}
v_j= \sqrt{ 1 - \frac{\bl_j^2}{n_j^2}} \;\; , \quad \bu_j = \frac{\bl_j}{n_j}
\ee
on the unit hemisphere $v^2 + \bu^2 =1\, ,\,v>0$ in $\cH_{n+1}$. If $|\bu_j|=1$ for some
$j$ we have another cusp in addition to the cusp at infinity, but now lying 
on the boundary $v=0$ of $\cH$. Therefore, the fundamental region always 
has the shape of a `skyscraper' that extends to infinite height
over the simplex $\Sigma\subset\R^n$ defined by the points $0$ and $\bu_j$, 
and whose `bottom' is cut off by the unit hemisphere.  See Fig.\,\ref{fig:lqs}
for an artist's view; the `bottom' of the skyscraper is the excised shaded region
on the unit sphere.

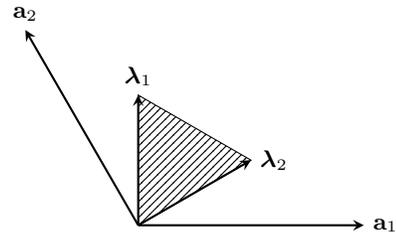
\begin{figure}
\begin{center}
\begin{tikzpicture}[>=stealth]

\filldraw[fill=gray!20!white, draw=black, pattern=north east lines] (0,0) -- (1.5,0.866) -- (0,1.732) -- (0,0);
\draw[->,thick] (0,0) -- (3,0) node[right]{$\ba_1$};
\draw[->,thick] (0,0) -- (-1.5,2.598) node[above]{$\ba_2$};
\draw[->,thick] (0,0) -- (0,1.732) node[above]{$\bl_1$};
\draw[->,thick] (0,0) -- (1.5,0.866) node[right]{$\bl_2$};

\end{tikzpicture}
\caption{Schematic example for the projection of the fundamental domain 
for the Weyl chamber onto the hypersurface $v=0$ for $A_2^{++}$. In
accordance with (\ref{lambda1}) the roots and weights here are normalized as 
$\ba_1=(1,0)$, $\ba_2=(-\frac{1}{2},\frac{\sqrt{3}}{2})$, and $\bl_1=(0,\frac{1}{\sqrt{3}})$, 
$\bl_2=(\frac{1}{2},\frac{1}{2\sqrt{3}})$.}\label{fig:projv}
\end{center}
\end{figure}

Using the above formulas we obtain
\be\label{S}
\bu_i\cdot \bu_j \equiv S_{ij} = \frac{1}{2m_i m_j} (B^{-1})_{ij} \ .
\ee
The matrix $S_{ij}$ encodes all the Lie algebraic 
information about the over-extended algebra $\mak g^{++}$ via the 
inverse symmetrized Cartan matrix $B^{-1}$ and the Coxeter labels $m_j$.
By a general result valid for all finite $\mak g$ \cite{Hum72} the matrices 
$B^{-1}$ are positive definite; furthermore their individual
entries $B^{-1}_{ij}$ are also positive.  It thus follows that
\be\label{Spos}
S > 0 \;\; \mbox{(as a matrix) and}\;\; S_{ij} > 0
\;\; \mbox{for all $i,j$}
\ee
Note that this formula holds for simply-laced as well as non-simply-laced 
(untwisted) algebras. In particular, in the non-simply laced case one has to 
distinguish between the Coxeter/dual Coxeter labels of the untwisted and 
the Coxeter/dual Coxeter labels of the twisted version of the over-extension of the algebra.

As we just explained the fundamental domain $\cF_0\subset\cH$ rises over the 
simplex $\Sigma  \subset\R^n$ defined by
\be
\Sigma :=
\big\{ \bX\in\R^n\, | \, \bX =\sum\limits_{i=1}^{n} t_i \bu_i\,;  \ t_i\geq 0,\ 
\sum\limits_{i=1}^{n} t_i\leq1 \big\} \ .
\ee
With the above definitions we get
\be
\bX(t)^2  =\sum\limits_{i,j=1}^n S_{ij}t_i t_j \ .
\ee
From the positivity properties (\ref{Spos})  we deduce the following 
chain of inequalities valid for all points $\bX(t)\in\Sigma$
\be\label{S>0}
0\leq \sum_{i,j} S_{ij} t_i t_j \,\leq\,  \max_{i,j} S_{ij} \left( \sum_k t_k\right)^2
       \leq \max_{i,j} S_{ij}
\ee
Therefore $\bX(t)^2 < 1$ as long as all matrix entries satisfy $S_{ij} < 1$. 
From \eqref{Lambda} it is straighforward to see that 
\be
\Lambda_j^2 = n_j^2 \,\big( \bu_j^2 -1 \big)
\ee
and it therefore follows that $S_{ii}=1$ when the corresponding 
hyperbolic weight $\Lambda_i$ becomes null; for spacelike weights
($\Lambda_j^2 >0$) we have $|\bu_j| > 1$, and the corresponding 
point $(v_j,\bu_j)$ no longer lies in the generalized upper half-plane.
This happens when $\mak g^{++}$ is Lorentzian, but no longer
hyperbolic, as is for instance the case for all $A_n^{++}$ with $n\geq 8$ 
and $B_n^{++}$ and $D_n^{++}$ for $n\geq 9$.

With the hyperbolic volume element
\be
\d{\rm vol} (\bu,v) = \frac{\d^nu\, \d
v}{v^{n+1}}
\ee
we thus obtain
\be
\vol(\mathcal{F}_0)=\sqrt{\det S}\int\limits_{\Delta_n}
\d t_1 \cdots \d t_n\int\limits_{\sqrt{1- \bX(t)^2}}^{\infty}\frac{\d v}{v^{n+1}}
\ee
where $\Delta_n$ is the {\em standard simplex} in $\R^n$ 
\be
\Delta_n := \big\{ (t_1,\dots,t_n) \, | \, t_i \geq 0 \, , \, \sum t_i \leq 1 \big\}
\ee
Performing the integral over $v$ we arrive at
\be\label{main}{}\nn\\
\boxed{\vol (\cF_0 ) = 
     \frac1{n}\int_{\Delta_n} \d t_1  \cdots \d t_n
\frac{\sqrt{\det S}}{(1-\sum t_i S_{ij}t_j)^{\frac{n}{2}}} \ .}\nn\\
\ee
This simple formula is our main result: it expresses the hyperbolic volume as an integral
over a standard simplex $\Delta_n$ in $\R^n$ with the single matrix $S_{ij}$ 
encoding all the Lie algebraic  information about the hyperbolic Weyl group. 
The integral is manifestly convergent if all $S_{ij} <1$. When $S_{ii}=1$
the corresponding point has $|\bu_i|=1$ and $v_i=0$ and thus lies on the 
boundary of $\cH$, but the integral is still convergent (see below for 
examples when this happens). For non-hyperbolic Lorentzian
algebras the integral diverges and therefore $\vol(\cF_0) = \infty$.

When evaluating this formula it may be convenient to diagonalize
the quadratic form in terms of new integration variables $\xi_i$
such that
\be
\sum_{i,j} S_{ij} t_i t_j = \xi_1^2 + \cdots + \xi_n^2
\ee
and the determinant factor $(\det S)^{1/2}$ is cancelled by the Jacobian.
The variables $\xi_i$ always exist by the positivity properties of the
matrix $S$. However, the (still simplicial) domain of integration is then 
more complicated to parametrize.

\section{Analytic Results}
We now show how our formula (\ref{main}) immediately yields the volumes
for various hyperbolic reflection groups corresponding to over-extended
hyperbolic algebras $\mak g^{++}$ of low rank. It is straightforward to check that
for $A_1$ (corresponding to the rank-3 Feingold-Frenkel algebra $A_1^{++}$)
we have $S=\frac12$, and one easily recovers the well known
result $\vol (\cF_0 [A_1^{++}])= \frac{\pi}6$. For this reason we proceed
right away to the case of rank 4.

The rank-4 algebras of over-extended type are $A_2^{++}, C_2^{++}$ and
$G_2^{++}$. For $\mak g^{++}= A_2^{++}$ we have $m_1 = m_2 =1$ and 
thus the matrix $S_{ij}$ is 1/2 the inverse of the $A_2$ Cartan matrix
\be
B^{-1}  =  \left( \begin{array}{cc} \frac23 & \frac13 \\[2mm]
                                   \frac13 & \frac23 \end{array}\right)
\;\; \Rightarrow \quad
S  =  \left( \begin{array}{cc} \frac13 & \frac16 \\[2mm]
                                   \frac16 & \frac13 \end{array}\right)
\ee
Transforming to new coordinates $\xi_1=\frac12(t_1 + t_2)$ and 
$\xi_2= (1/2\sqrt{3}) (t_1 -t_2)$ such that the Jacobi determinant
cancels the factor $(\det S)^{1/2}= 1/2\sqrt{3}$ and
\be\label{b2} 
\frac13 \big( t_1^2 + t_1 t_2 + t_2^2 \big) = \xi_1^2 + \xi_2^2
\ee
we obtain
\be
\vol (\cF_0[A_2^{++}]) &=& \frac12 
\int_{0}^\frac{1}{2} \d \xi_1 \int_{- \frac{\xi_1}{\sqrt{3}}}^{\frac{\xi_1}{\sqrt{3}}}
\frac{\d\xi_2}{1- \xi_1^2 - \xi_2^2} \nn\\
&&\hspace{-2cm}
 =  \frac{1}{2}\int_0^\frac{1}{2}\frac{\d \xi}{\sqrt{1-\xi^2}}\ln \left( 
     \frac{\sqrt{1-\xi^2}+\frac{1}{\sqrt{3}}\xi}{\sqrt{1-\xi^2} -\frac{1}{\sqrt{3}}\xi} \right)
\ee
The substitution $ \xi =\sin\theta$ leads to
\begin{align}
\vol(\cF_0[A_2^{++}]) & = 
\frac{1}{2} \int_{0}^{\frac{\pi}{6}}\d\theta \ln \left(\frac{\cos\theta+\frac{1}{\sqrt{3}}\sin\theta}{\cos\theta-\frac{1}{\sqrt{3}}\sin\theta}\right) \nn\\
& \hspace{-1cm}= \frac{1}{2} \int_{0}^{\frac{\pi}{6}}
    \d\theta \ln \left(\frac{2\sin(\theta+\frac{\pi}{3})}{2\sin(\frac{\pi}{3} - \theta)}\right) 
\end{align}
After a suitable shift of integration variables and using
the definition and properties of the Lobachevsky function, 
this reduces to
\selectlanguage{russian}
\be
& \phantom{=} & \hskip-3em \vol(\cF_0[A_2^{++}])\\
&=&\frac12 \left[ \Lo\left(\frac{\pi}3\right)  - \Lo\left(\frac{\pi}6\right)
     - \Lo\left(\frac{\pi}2\right)  + \Lo\left(\frac{\pi}3\right)\right] \nn\\
&=&   \frac14   \,  \Lo\left(\frac{\pi}3\right)
\ee

For $\mathfrak{g}^{++} = G_2^{++}$ we have the Coxeter labels $m_1 =2 \, , \, m_2 = 3$ 
and the relevant matrices are
\be
B^{-1}  = \left( \begin{array}{cc} 2 & 3 \\[2mm]
                                                               3 & 6 \end{array}\right)
\;\; \Rightarrow \quad
S  =   \left( \begin{array}{cc} \frac14 & \frac14 \\[2mm]
                                   \frac14 & \frac13 \end{array}\right)                
\ee
Now the substitution to diagonalize the quadratic form is
$\xi_1=\frac12(t_1 + t_2)\, , \, \xi_2= (1/2\sqrt{3}) t_2$, and we get
\be
\vol (\cF_0[G_2^{++}]) &=&  \frac12
\int_{0}^\frac{1}{2} \d \xi_1 \int_0^{\frac{\xi_1}{\sqrt{3}}}
\frac{\d \xi_2}{1- \xi_1^2 - \xi_2^2} \nn\\
 && \hspace{-1.5cm} =  \frac12 \, \vol \big( \cF_0 [A_2^{++}] \big)
 = \frac18 \, \Lo\left(\frac{\pi}3\right)
\ee
Finally, for $C_2^{++}$ we have $m_1=m_2 = 1$ and
\be
B^{-1} = \left( \begin{array}{cc} \frac12 & \frac12 \\[2mm]
                                                        \frac12 & 1 \end{array}\right)
\;\; \Rightarrow \quad
S  =  \left( \begin{array}{cc} \frac14& \frac14 \\[2mm]
                                   \frac14 & \frac12 \end{array}\right)
\ee
Note that the corresponding matrix for $B_2^{++}$ is
\mbox{$S  =  \frac14 \left(\begin{smallmatrix} 2& 1 \\ 1 & 1 \end{smallmatrix}\right)$}
and yields the same volume.
\\
Now we substitute $\xi_1 = \frac12(t_1 + t_2)$ and $\xi_2 = \frac12 t_2$ to get
\be
\vol \big (\cF_0[C_2^{++}] \big) &=&  \frac12 
\int_{0}^\frac{1}{2} \d \xi_1 \int_0^{\xi_1} \frac{\d\xi_2}{1- \xi_1^2- \xi_2^2} \nn\\ 
&& \hspace{-1.5cm}= \frac{1}{2} \int_{0}^{\frac{\pi}{6}} 
    \d\theta \ln \left(\frac{2\sin(\theta+\frac{\pi}{4})}{2\sin(\frac{\pi}{4} - \theta)}\right) 
\ee
Similar manipulations as before lead to the result
\be
& \phantom{=} & \hskip-3em \vol\left( \cF_0[C_2^{++}] \right) \nn\\ 
& = & \frac14 \left[ 
\Lo\left(\frac{\pi}4\right) - \Lo\left(\frac{\pi}{12}\right) - \Lo\left(\frac{5\pi}{12}\right) +  \Lo\left(\frac{\pi}4\right) \right] \nn\\ 
&=&\frac16 \, \Lo\left(\frac{\pi}{4}\right)
\ee

\section{Higher rank algebras}

What about higher rank algebras? Although the integrals (\ref{main}) look 
elementary it turns out that calculations become rapidly more complicated with 
increasing dimension, and we have not been able to derive `simple' closed form expressions
for them when $n>2$. The complications are mainly due to the integration boundaries which must 
be analyzed case by case.  Although (\ref{main}) is suggestive of higher order Lobachevsky functions (see appendix), this expectation (as expressed, for instance, in
\cite{Vin93b}) is not borne out by the concrete calculations, nor have such expressions 
been explicitly exhibited in the literature, see e.g. \cite{Vin93b}. One possibility, to
be explored in future work, would be to expand the integrand in (\ref{main}) whereby
the integral is expressed as an infinite sum of terms each one of which involves
an integral of a monomial over the standard simplex $\Delta_n$. Such integrals have been studied in the literature \cite{Bri88,BBDLKV11} but the resulting expressions are still quite 
involved. Numerically these series converge rapidly, as all terms are of the same sign.

Using \eqref{main} one can compute the volume of different fundamental domains 
numerically. The only input information that is needed is the matrix $S$, which is 
calculated from the matrix $B^{-1}$ and the Coxeter labels via \eqref{S}. As already 
mentioned above, $B^{-1}$ is the inverse symmetrized Cartan matrix and its form for the different algebras can be found in the standard Lie algebra literature, see e.g. \cite{Hum72}).
\selectlanguage{english}
\begin{longtable*}[b]{C{1cm}|C{7cm}|C{7cm}}
$G$ & Untwisted & Twisted \\ \hline\hline
\begin{minipage}[b]{1cm}
\centering $B_n$\\ \phantom{ }
\end{minipage}
&
\begin{tikzpicture}
[place/.style={circle,draw=black,fill=black,
inner sep=0pt,minimum size=6}]
\draw (0,1) -- (0,0) -- (1,0);
\draw[dashed] (1,0) -- (2,0);
\draw (2,0.1) -- (3,0.1);
\draw (2,-.1) -- (3,-.1);
\draw (-2,0) -- (0,0);
\draw (2.6,0) -- (2.4,.2);
\draw (2.6,0) -- (2.4,-.2);
\node at (-2,0) [place,label=below:$-1$] {};
\node at (0,1) [place,label=right:$1$] {};
\node at (-1,0) [place,label=below:$0$] {};
\node at (0,0) [place,label=below:$2$] {};
\node at (1,0) [place,label=below:$3$] {};
\node at (2,0) [place,label=below:$n-1$] {};
\node at (3,0) [place,label=below:$n\!\!\!\phantom{1}$] {};
\end{tikzpicture}
&
\begin{tikzpicture}
[place/.style={circle,draw=black,fill=black,
inner sep=0pt,minimum size=6}]
\draw (-1,0) -- (0,0);
\draw (0,.1) -- (1,.1);
\draw (0,-.1) -- (1,-.1);
\draw (1,0) -- (2,0);
\draw[dashed] (2,0) -- (3,0);
\draw (3,.1) -- (4,.1);
\draw (3,-.1) -- (4,-.1);
\draw (0.4,0) -- (0.6,0.2);
\draw (0.4,0) -- (0.6,-0.2);
\draw (3.6,0) -- (3.4,0.2);
\draw (3.6,0) -- (3.4,-0.2);
\node at (-1,0) [place,label=below:$-1$] {};
\node at (0,0) [place,label=below:$0$] {};
\node at (1,0) [place,label=below:$1$] {};
\node at (2,0) [place,label=below:$2$] {};
\node at (3,0) [place,label=below:$n-1$] {};
\node at (4,0) [place,label=below:$n\!\!\!\phantom{1}$] {};
\end{tikzpicture}
\\ \hline
\begin{minipage}[b]{1cm}
\centering $C_n$\\ \phantom{ }
\end{minipage} & \begin{tikzpicture}
[place/.style={circle,draw=black,fill=black,
inner sep=0pt,minimum size=6}]
\draw (-1,0) -- (0,0);
\draw (0,.1) -- (1,.1);
\draw (0,-.1) -- (1,-.1);
\draw (1,0) -- (2,0);
\draw[dashed] (2,0) -- (3,0);
\draw (3,.1) -- (4,.1);
\draw (3,-.1) -- (4,-.1);
\draw (0.6,0) -- (0.4,0.2);
\draw (0.6,0) -- (0.4,-0.2);
\draw (3.4,0) -- (3.6,0.2);
\draw (3.4,0) -- (3.6,-0.2);
\node at (-1,0) [place,label=below:$-1$] {};
\node at (0,0) [place,label=below:$0$] {};
\node at (1,0) [place,label=below:$1$] {};
\node at (2,0) [place,label=below:$2$] {};
\node at (3,0) [place,label=below:$n-1$] {};
\node at (4,0) [place,label=below:$n\!\!\!\phantom{1}$] {};
\end{tikzpicture}
&
\begin{tikzpicture}
[place/.style={circle,draw=black,fill=black,
inner sep=0pt,minimum size=6}]
\draw (0,1) -- (0,0) -- (1,0);
\draw[dashed] (1,0) -- (2,0);
\draw (2,0.1) -- (3,0.1);
\draw (2,-.1) -- (3,-.1);
\draw (-2,0) -- (0,0);
\draw (2.4,0) -- (2.6,.2);
\draw (2.4,0) -- (2.6,-.2);
\node at (-2,0) [place,label=below:$-1$] {};
\node at (0,1) [place,label=right:$1$] {};
\node at (-1,0) [place,label=below:$0$] {};
\node at (0,0) [place,label=below:$2$] {};
\node at (1,0) [place,label=below:$3$] {};
\node at (2,0) [place,label=below:$n-1$] {};
\node at (3,0) [place,label=below:$n\!\!\!\phantom{1}$] {};
\end{tikzpicture} 
\\ \hline
\rule{0pt}{1.2cm}
\begin{minipage}[b]{1cm}
\centering $F_4$\\ \phantom{ }
\end{minipage} & \begin{tikzpicture}
[place/.style={circle,draw=black,fill=black,
inner sep=0pt,minimum size=6}]
\draw (-1,0) -- (2,0);
\draw (2,.1) -- (3,.1);
\draw (2,-.1) -- (3,-.1);
\draw (3,0) -- (4,0);
\draw (2.6,0) -- (2.4,0.2);
\draw (2.6,0) -- (2.4,-0.2);
\node at (-1,0) [place,label=below:$-1$] {};
\node at (0,0) [place,label=below:$0$] {};
\node at (1,0) [place,label=below:$1$] {};
\node at (2,0) [place,label=below:$2$] {};
\node at (3,0) [place,label=below:$3$] {};
\node at (4,0) [place,label=below:$4$] {};
\end{tikzpicture}
&
\begin{tikzpicture}
[place/.style={circle,draw=black,fill=black,
inner sep=0pt,minimum size=6}]
\draw (-1,0) -- (2,0);
\draw (2,.1) -- (3,.1);
\draw (2,-.1) -- (3,-.1);
\draw (3,0) -- (4,0);
\draw (2.4,0) -- (2.6,0.2);
\draw (2.4,0) -- (2.6,-0.2);
\node at (-1,0) [place,label=below:$-1$] {};
\node at (0,0) [place,label=below:$0$] {};
\node at (1,0) [place,label=below:$1$] {};
\node at (2,0) [place,label=below:$2$] {};
\node at (3,0) [place,label=below:$3$] {};
\node at (4,0) [place,label=below:$4$] {};
\end{tikzpicture}
\\ \hline
\rule{0pt}{1.2cm}
\begin{minipage}[b]{1cm}
\centering $G_2$\\ \phantom{ }
\end{minipage} & \begin{tikzpicture}
[place/.style={circle,draw=black,fill=black,
inner sep=0pt,minimum size=6}]
\draw (-1,0) -- (2,0);
\draw (1,.1) -- (2,.1);
\draw (1,-.1) -- (2,-.1);
\draw (1.6,0) -- (1.4,0.2);
\draw (1.6,0) -- (1.4,-0.2);
\node at (-1,0) [place,label=below:$-1$] {};
\node at (0,0) [place,label=below:$0$] {};
\node at (1,0) [place,label=below:$1$] {};
\node at (2,0) [place,label=below:$2$] {};
\end{tikzpicture}
&
\begin{tikzpicture}
[place/.style={circle,draw=black,fill=black,
inner sep=0pt,minimum size=6}]
\draw (-1,0) -- (2,0);
\draw (1,.1) -- (2,.1);
\draw (1,-.1) -- (2,-.1);
\draw (1.4,0) -- (1.6,0.2);
\draw (1.4,0) -- (1.6,-0.2);
\node at (-1,0) [place,label=below:$-1$] {};
\node at (0,0) [place,label=below:$0$] {};
\node at (1,0) [place,label=below:$1$] {};
\node at (2,0) [place,label=below:$2$] {};
\end{tikzpicture}\\ \hline\hline
\multicolumn{3}{c}{} \\[1mm]
\caption{Dynkin diagrams of the over-extended twisted and untwisted non-simply laced finite-dimensional algebras with Dynkin labeling of nodes}\label{fig:nsla}
\end{longtable*}
\selectlanguage{russian}
Here we list the matrices $S$ for the Lie Algebras of 
$A_n, B_n, C_n, D_n, F_4, E_6, E_7$ 
and $E_8$ in the Cartan classification (the matrices for the rank 2 algebras were already
given in the previous section). In addition we list the set of Coxeter (or dual Coxeter) 
labels $m_i$ used in the computation of $S$. Note that in the case of the non-simply laced algebras it is necessary to distinguish between the labels of the \textit{twisted} and \textit{untwisted} algebra. For these algebras we label the matrix $S$ with a superscript $^{(1)}$ 
or $^{(2)}$, respectively, indicating whether it corresponds to the \textit{untwisted} or 
\textit{twisted} over-extension of the algebra. Considering Table \ref{fig:nsla} containing Dynkin diagrams of over-extensions of the non-simply laced algebras, we note that the twisted Dynkin diagram of $B_n$ is the same as the untwisted Dynkin diagram of $C_n$, simply with the direction of the arrows reversed. This tells us that the volumes of the corresponding fundamental domains have to be the same, since the matrix $S$ is the \textit{symmetrized} version the Cartan matrix and therefore contains no information about the direction of the arrows. A similar correspondence holds for the untwisted diagram of $B_n$ and the twisted diagram of $C_n$, as well as for $G_2$ and $F_4$.\\
\\
The condition for the over-extension of each algebra to be of hyperbolic type is that \textit{all} of the diagonal entries $S_{ii}$ of the underlying finite-dimensional algebra must satisfy $S_{ii}\leq1$. Geometrically each $S_{ii}=1$ corresponds to an additional fundamental weight (edge of the Weyl chamber) lying \textit{on} the forward light-cone. For each $S_{ii}>1$ an additional fundamental weight lies \textit{outside} the light cone and the over-extension is not of hyperbolic type. For each algebra we state the range of $n$ for which the over-extension is hyperbolic.\\
\\
\noindent{$\mak{g}^{++} \! = \!A_n^{++}$:}
The Coxeter labels are $m_i=\left(1,...,1\right)$, and thus
the matrix $S$ is
\begin{align}
& S[A_n]=\frac{1}{(n+1)}\times \nn\\
&\times\begin{pmatrix}
\frac{n}{2} & \frac{n-1}{2} & \frac{n-2}{2} & \frac{n-3}{2} &  \cdots & \frac12 \\[2mm]
\frac{n-1}{2} & n-1 & n-2 & n-3 & \cdots & 1 \\[2mm]
\frac{n-2}{2} & n-2 & \frac{3(n-2)}{2} & \frac{3(n-3)}{2} & \cdots & \frac32 \\[2mm]
\frac{n-3}{2} & n-3 & \frac{3(n-3)}{2} & 2(n-3) & \cdots & 2 \\[2mm]
\vdots & \vdots & \vdots & \vdots & \ddots & \vdots\\[2mm]
\frac12 & 1 & \frac32 & 2 & \cdots & \frac{n}{2}	
\end{pmatrix}
\end{align}
From the explicit form of the matrix it is obvious that 
\be
S_{ij}\leq 1 \quad \Leftrightarrow\quad 
\frac{j(n+1-j)}{2(n+1)} \leq 1
\ee
for all $j=1,\dots, n$. Hence $A_n^{++}$ is hyperbolic for $n\leq 7$.

\vspace{1mm}

\noindent $\mak{g}^{++} \!=\! B_n^{++}$: the Coxeter labels 
are $m_i=\left(1,2,...,2\right)$ (as untwisted Coxeter labels for $B_n$, 
and as twisted dual Coxeter labels for $C_n$), and the matrix $S$ is
\begin{align}
S^{(1)} [B_n] & =\left( \begin{array}{c|cccccc}
\frac12 & \frac14 & \frac14 & \frac14 & \cdots & \frac14 & \frac14 \\[2mm] \mhlines
\frac14 & \frac14 & \frac14 & \frac14 & \cdots & \frac14 & \frac14 \\[2mm]
\frac14 & \frac14 & \frac38 & \frac38 & \cdots & \frac38 & \frac38 \\[2mm]
\frac14 & \frac14 & \frac38 & \frac12 & \cdots & \frac12 & \frac12 \\[2mm]
\vdots & \vdots & \vdots & \vdots & \ddots & \vdots & \vdots \\[2mm]
\frac14 & \frac14 & \frac38 & \frac12 & \cdots & \frac{n-1}{8} & \frac{n-1}{8} \\[2mm]
\frac14 & \frac14 & \frac38 & \frac12 & \cdots & \frac{n-1}{8} & \frac{n}{8}
\end{array}\right)
\end{align}
All matrix entries are $\leq 1$ for $n\leq 8$, whence
$B_n^{++}$ is hyperbolic for $n\leq8$. Inverting the arrow
in the Dynkin diagram we infer that
\be
S^{(1)}[B_n] = S^{(2)}[C_n].
\ee

\vspace{1mm}

\noindent $\mak{g}^{++} \!=\! C_n^{++}$: the Coxeter labels are
$m_i=\left(1,...,1\right)$ (as twisted Coxeter labels for $B_n$ and 
untwisted dual Coxeter labels for $C_n$), so
\begin{align}
S^{(1)} [C_n] & =\begin{pmatrix}
\frac14 & \frac14 & \frac14 & \frac14 & \cdots & \frac14 & \frac14 \\[2mm]
\frac14 & \frac12 & \frac12 & \frac12 & \cdots & \frac12 & \frac12 \\[2mm]
\frac14 & \frac12 & \frac34 & \frac34 & \cdots & \frac34 & \frac34 \\[2mm]
\frac14 & \frac12 & \frac34 & 1 & \cdots & 1 & 1 \\[2mm]
\vdots & \vdots & \vdots & \vdots & \ddots & \vdots & \vdots \\[2mm]
\frac14 & \frac12 & \frac34 & 1 & \cdots & \frac{n-1}{4} & \frac{n-1}{4} \\[2mm]
\frac14 & \frac12 & \frac34 & 1 & \cdots & \frac{n-1}{4} & \frac{n}{4}
\end{pmatrix}
\end{align}
Clearly,  $C_n^{++}$ is hyperbolic for $n\leq4$.  As before we get
\be
S^{(1)}[C_n] = S^{(2)}[B_n]
\ee

\vspace{1mm}

\noindent $\mak{g}^{++} \!=\! D_n^{++}$: the Coxeter labels are
$m_i=\left(1,2,...,2,1,1\right)$, and therefore
\be
S[D_n]  = \left(\begin{array}{c|ccccc|cc}
\frac12 & \frac14 & \frac14 & \frac14 & \cdots & \frac14 & \frac14 & \frac14 \\[2mm] \mhlinee
\frac14 & \frac14 & \frac14 & \frac14 & \cdots & \frac14 & \frac14 & \frac14 \\[2mm]
\frac14 & \frac14 & \frac38 & \frac38 & \cdots & \frac38 & \frac38 & \frac38 \\[2mm]
\frac14 & \frac14 & \frac38 & \frac12 & \cdots & \frac12 & \frac12 & \frac12 \\[2mm]
\vdots & \vdots & \vdots & \vdots & \ddots & \vdots & \vdots  & \vdots \\[2mm]
\frac14 & \frac14 & \frac38 & \frac12 & \cdots & \frac{n-2}{8} & \frac{n-2}{8} & \frac{n-2}{8} \\[2mm]\mhlinee
\frac14 & \frac14 & \frac38 & \frac12 & \cdots & \frac{n-2}{8} & \frac{n}{8} & \frac{n-2}{8} \\[2mm]
\frac14 & \frac14 & \frac38 & \frac12 & \cdots & \frac{n-2}{8} & \frac{n-2}{8} & \frac{n}{8}
\end{array}\right)
\ee
We see that $D_n^{++}$ is hyperbolic for $n\leq 8$.

\vspace{2mm}

\noindent $\mak{g}^{++} \!=\! F_4^{++}$: the Coxeter labels are
$\left(2,3,2,1\right)$ for the untwisted dual Coxeter labels as 
well as for the twisted Coxeter labels:
\be
S^{(1)}[F_4] = S^{(2)}[F_4] = \left(\begin{array}{cccc} \frac14 & \frac14 & \frac14 & \frac14 \\[2mm]
						  \frac14 & \frac13 & \frac13 & \frac13 \\[2mm]
						  \frac14 & \frac13 & \frac38 & \frac38 \\[2mm]
						  \frac14 & \frac13 & \frac38 & \frac12 \end{array} \right) 
\ee

\vspace{2mm}

\noindent $\mak{g}^{++} \!=\! E_6^{++}$: the Coxeter labels are
$\left(1,2,3,2,1,2\right)$, and we have
\be
S[E_6] = \left(\begin{array}{cccccc} \frac23 & \frac{5}{12} & \frac13 & \frac13 & \frac13 & \frac14 \\[2mm]
						\frac{5}{12} & \frac{5}{12} & \frac13 & \frac13 & \frac13 & \frac14 \\[2mm]
						\frac13 & \frac13 & \frac13 & \frac13 & \frac13 & \frac14 \\[2mm]
						\frac13 & \frac13 & \frac13 & \frac{5}{12} & \frac{5}{12} & \frac14 \\[2mm]
						\frac13 & \frac13 & \frac13 & \frac{5}{12} & \frac23 & \frac14 \\[2mm]
						\frac14 & \frac14 & \frac14 & \frac14 & \frac14 & \frac14 \end{array} \right) 
\ee

\vspace{2mm}

\noindent $\mak{g}^{++} \!=\! E_7^{++}$: the Coxeter labels are
$\left(2,3,4,3,2,1\right)$, and we have
\be
S[E_7] = \left(\begin{array}{ccccccc} \frac14 & \frac14 & \frac14 & \frac14 & \frac14 & \frac14 & \frac14 \\[2mm]
						\frac14 & \frac13 & \frac13 & \frac13 & \frac13 & \frac13 & \frac13 \\[2mm]
						\frac14 & \frac13 & \frac38 & \frac38 & \frac38 & \frac38 & \frac38 \\[2mm]
						\frac14 & \frac13 & \frac38 & \frac{5}{12} & \frac{5}{12} & \frac{5}{12} & \frac38 \\[2mm]
						\frac14 & \frac13 & \frac38 & \frac{5}{12} & \frac12 & \frac12 & \frac38 \\[2mm]
						\frac14 & \frac13 & \frac38 & \frac{5}{12} & \frac12 & \frac34 & \frac38 \\[2mm] 
						\frac14 & \frac13 & \frac38 & \frac38 & \frac38 & \frac38 & \frac{7}{16} \end{array} \right) 
\ee
\\

\vspace{2mm}

\noindent $\mak{g}^{++} \!=\! E_8^{++}$: with the Coxeter labels 
$\left(2,3,4,5,6,4,2,3\right)$ we have
\be
S[E_8] =\left(\begin{array}{cccccccc} \frac14 & \frac14 & \frac14 & \frac14 & \frac14 & \frac14 & \frac14 & \frac14 \\[2mm]
						\frac14 & \frac13 & \frac13 & \frac13 & \frac13 & \frac13 & \frac13 & \frac13 \\[2mm]
						\frac14 & \frac13 & \frac38 & \frac38 & \frac38 & \frac38 & \frac38 & \frac38 \\[2mm]
						\frac14 & \frac13 & \frac38 & \frac25 & \frac25 & \frac25 & \frac25 & \frac25 \\[2mm]
						\frac14 & \frac13 & \frac38 & \frac25 & \frac{5}{12} & \frac{5}{12} & \frac{5}{12} & \frac{5}{12} \\[2mm]
						\frac14 & \frac13 & \frac38 & \frac25 & \frac{5}{12} & \frac{7}{16} & \frac{7}{16} & \frac{5}{12} \\[2mm] 
						\frac14 & \frac13 & \frac38 & \frac25 & \frac{5}{12} & \frac{7}{16} & \frac12 & \frac{5}{12} \\[2mm]
						\frac14 & \frac13 & \frac38 & \frac25 & \frac{5}{12} & \frac{5}{12} & \frac{5}{12} & \frac29  \end{array} \right) 
\ee

\vspace{2mm}

Using equations \eqref{UHPMap} one can determine the coordinates of 
the vertices of the fundamental domain. For example, the vertices of the domain 
corresponding to the Weyl group of $\mak{e}_{10}$ are given by
\be
\left(v_1,\bu_1\right)& = &\left(\frac{\sqrt{3}}{2},\frac12 \beb_0 \right) \nonumber \\
\left(v_2,\bu_2\right)& = &\left(\sqrt{\frac23},\frac12 \beb_0+\frac16\left(\beb_1+\beb_5+\beb_6\right)\right) \nonumber \\
\left(v_3,\bu_3\right)& = &\left(\sqrt{\frac58},\frac12 \beb_0+\frac14\left(\beb_5+\beb_6\right)\right) \nonumber \\
\left(v_4,\bu_4\right)& = &\left(\sqrt{\frac35},\frac12 \beb_0+\frac{1}{10}\left(\beb_2+3\beb_5+2\beb_6-\beb_7\right) \right)\nonumber \\
\left(v_5,\bu_5\right)& = &\left(\frac{1}{\sqrt{6}},\frac12 \beb_0+\frac16\left(2\beb_5+\beb_6-\beb_7\right) \right) \nonumber \\
\left(v_6,\bu_6\right)& = &\left(\frac34,\frac12 \beb_0+\frac18\left(\beb_3+3\beb_5+\beb_6-\beb_7\right) \right) \nonumber \\
\left(v_7,\bu_7\right)& = &\left(\frac{1}{\sqrt{2}},\frac12 \beb_0+\frac12\beb_5 \right) \nonumber \\
\left(v_8,\bu_8\right)& = &\left(\frac{\sqrt{5}}{3},\frac12 \beb_0+\frac16\left(\beb_4+2\beb_5+\beb_6-\beb_7\right) \right)
\ee
The special feature of this example is that the vectors $\bu_j$ now belong to
the octonions $\oct$, the non-commutative and non-associative maximal division
algebra. Accordingly, the unit vectors $\beb_j$ (for $j=1,\dots,7$)
are just the octonionic imaginary units. The vertex coordinates of the fundamental 
domains of other Weyl groups are obtained similarly.

By evaluating the integrals in \eqref{main} numerically we obtain the volumes of all the fundamental domains of the hyperbolic Weyl groups of the algebras listed above. Employing a deterministic adaptive integration scheme with a sufficient number of evaluation points of the integrand, the 
values we find agree to high accuracy with those found in \cite{JKRT99} where the 
volumes of all hyperbolic Coxeter simplices were obtained by a different method.

\acknowledgments{We are very grateful to Axel Kleinschmidt for discussions
and helpful comments on an earlier version of this paper. We would also like
to thank Jakob Palmkvist  for discussions and correspondence. The work of
P. Fleig is supported by an IRAP Erasmus Mundus Joint Doctorate Fellowship 
and the University of Nice -- Sophia Antipolis.}

\vspace{1cm}

\appendix

\selectlanguage{english}
\section{The Lobachevsky function}
\selectlanguage{russian}

The Lobachevsky function $\Lo$ is defined as
\be
\Lo (\theta) = -\int_0^\theta \log (|2\sin t |)\d t \quad \forall\theta\in\R
\ee
and related to the dilogarithm and the Clausen function through the identities
\be\label{lilo}
\Lo (\omega ) = \frac{1}{2} \Im (\Li_2(\e^{2\I\omega}))=\frac{1}{2}\mathrm{Cl}_2(2\omega) \ .
\ee
where the dilogarithm is defined as
\be
\Li_2(z)=-\int_0^z\frac{\log(1-w)}{w}\d w \quad \forall z\in\C: z\notin [1,\infty) \ .
\ee
The Lobachevsky function satisfies the relations
\be
\Lo(0)=\Lo\left(\frac{\pi}{2}\right)=0 \\
\Lo(\theta+\pi)=\Lo(\theta) \\
\Lo(-\theta)=-\Lo(\theta)
\ee
and
\be
\Lo(n\theta)=n\sum_{j=0}^{n-1}\Lo(\theta+\frac{j\pi}{n})\quad \forall n\in\Z_+ \ ,
\ee
yielding e.g.
\be
\Lo\left(\frac{\pi}{6}\right) &=&\frac{3}{2}\Lo\left(\frac{\pi}{3}\right)  \nn\\
\Lo\left(\frac{\pi}{4}\right) &=& \frac34\left[ 
\Lo\left(\frac{\pi}{12}\right) \ + \Lo\left(\frac{5\pi}{12}\right) \right]\ .
\ee
The polylogarithm functions
\be
\Li_n (z) = \sum_{r=1}^{\infty}\frac{z^r}{r^n}  \quad\forall z\in\C\;\;\;\;\forall n\geq 1
\ee
arise naturally in the computation of hyperbolic volume. They are inductively related through
\be
\Li_1 (z) = -\log (1-z) \\
\Li_n (z) =\int_0^z \Li_{n-1} (w) \frac{\d w}{w} \ .
\ee
It is furthermore common to define the higher Lobachevsky functions through the polylogarithm according to
\be
\Lo_{2m}(\theta)=\frac{1}{2^{2m-1}}\Im(\Li_{2m}(\e^{2\I\theta}))\\
\Lo_{2m+1}(\theta)=\frac{1}{2^{2m}}\Re(\Li_{2m+1}(\e^{2\I\theta})) \ .
\ee
The higher Lobachevsky functions satisfy the more general relations
\be
\Lo_m(\theta)=\Lo_m(\theta+\pi)\\
\frac{1}{n^{m-1}}\Lo_m(n\theta)= \sum_{j=0}^{n-1}\Lo_m\left(\theta+\frac{j\pi}{n}\right) \\
\Lo_m(-\theta)=(-1)^{m+1}\Lo_m(\theta) \ .
\ee
However, as we already pointed out, and unlike for rank four, it does not appear that 
the volumes for the higher rank fundamental domains can be expressed solely in
terms of higher Lobachevsky functions.

\bibliography{volume}

\begin{thebibliography}{16}
\expandafter\ifx\csname natexlab\endcsname\relax\def\natexlab#1{#1}\fi
\expandafter\ifx\csname bibnamefont\endcsname\relax
  \def\bibnamefont#1{#1}\fi
\expandafter\ifx\csname bibfnamefont\endcsname\relax
  \def\bibfnamefont#1{#1}\fi
\expandafter\ifx\csname citenamefont\endcsname\relax
  \def\citenamefont#1{#1}\fi
\expandafter\ifx\csname url\endcsname\relax
  \def\url#1{\texttt{#1}}\fi
\expandafter\ifx\csname urlprefix\endcsname\relax\def\urlprefix{URL }\fi
\providecommand{\bibinfo}[2]{#2}
\providecommand{\eprint}[2][]{\url{#2}}

\bibitem[{\citenamefont{Apostol}(1990)}]{Apo90}
\bibinfo{author}{\bibfnamefont{T.~M.} \bibnamefont{Apostol}},
  \emph{\bibinfo{title}{Modular Functions and Dirichlet Series in Number Theory
  (2nd ed)}}, vol.~\bibinfo{volume}{41} of \emph{\bibinfo{series}{Graduate
  Texts in Mathematics}} (\bibinfo{publisher}{Springer-Verlag, New York},
  \bibinfo{year}{1990}).

\bibitem[{\citenamefont{Kac}(1990)}]{Kac90}
\bibinfo{author}{\bibfnamefont{V.~G.} \bibnamefont{Kac}},
  \emph{\bibinfo{title}{Infinite Dimensional Lie Algebras (3rd ed.)}}
  (\bibinfo{publisher}{Cambridge University Press}, \bibinfo{year}{1990}).

\bibitem[{\citenamefont{Feingold et~al.}(2009)\citenamefont{Feingold,
  Kleinschmidt, and Nicolai}}]{FKN09}
\bibinfo{author}{\bibfnamefont{A.~J.} \bibnamefont{Feingold}},
  \bibinfo{author}{\bibfnamefont{A.}~\bibnamefont{Kleinschmidt}},
  \bibnamefont{and} \bibinfo{author}{\bibfnamefont{H.}~\bibnamefont{Nicolai}},
  \bibinfo{journal}{J Algebra} \textbf{\bibinfo{volume}{322}},
  \bibinfo{pages}{1295} (\bibinfo{year}{2009}), \eprint{arXiv:0805.3018
  [math.RT]}.

\bibitem[{\citenamefont{Feingold and Frenkel}(1983)}]{FF}
\bibinfo{author}{\bibfnamefont{A.~J.} \bibnamefont{Feingold}} \bibnamefont{and}
  \bibinfo{author}{\bibfnamefont{I.~B.} \bibnamefont{Frenkel}},
  \bibinfo{journal}{Math Ann} \textbf{\bibinfo{volume}{263}},
  \bibinfo{pages}{87} (\bibinfo{year}{1983}).

\bibitem[{\citenamefont{Kleinschmidt et~al.}(2010)\citenamefont{Kleinschmidt,
  Nicolai, and Palmkvist}}]{KNP10}
\bibinfo{author}{\bibfnamefont{A.}~\bibnamefont{Kleinschmidt}},
  \bibinfo{author}{\bibfnamefont{H.}~\bibnamefont{Nicolai}}, \bibnamefont{and}
  \bibinfo{author}{\bibfnamefont{J.}~\bibnamefont{Palmkvist}}
  (\bibinfo{year}{2010}), \eprint{arXiv:1010.2212v1 [math.NT]}.

\bibitem[{\citenamefont{Lobachevsky}(1904)}]{Lob04}
\bibinfo{author}{\bibfnamefont{N.~I.} \bibnamefont{Lobachevsky}},
  \bibinfo{journal}{Deutsche \"Ubersetzung von H. Liebmann, Teubner, Leipzig}
  (\bibinfo{year}{1904}).

\bibitem[{\citenamefont{Coxeter}(1935)}]{Cox35}
\bibinfo{author}{\bibfnamefont{H.}~\bibnamefont{Coxeter}}, \bibinfo{journal}{Q
  J Math} \textbf{\bibinfo{volume}{6}}, \bibinfo{pages}{13}
  (\bibinfo{year}{1935}).

\bibitem[{\citenamefont{Vinberg}(1993)}]{Vin93b}
\bibinfo{editor}{\bibfnamefont{E.~B.} \bibnamefont{Vinberg}}, ed.,
  \emph{\bibinfo{title}{Geometry II: Spaces of Constant Curvature}}
  (\bibinfo{publisher}{Springer-Verlag}, \bibinfo{year}{1993}).

\bibitem[{\citenamefont{Johnson et~al.}(1999)\citenamefont{Johnson, Kellerhals,
  Ratcliffe, and Tschantz}}]{JKRT99}
\bibinfo{author}{\bibfnamefont{N.}~\bibnamefont{Johnson}},
  \bibinfo{author}{\bibfnamefont{R.}~\bibnamefont{Kellerhals}},
  \bibinfo{author}{\bibfnamefont{J.}~\bibnamefont{Ratcliffe}},
  \bibnamefont{and} \bibinfo{author}{\bibfnamefont{S.}~\bibnamefont{Tschantz}},
  \bibinfo{journal}{Transform Groups} \textbf{\bibinfo{volume}{4}},
  \bibinfo{pages}{329} (\bibinfo{year}{1999}).

\bibitem[{\citenamefont{Kleinschmidt et~al.}(2009)\citenamefont{Kleinschmidt,
  Koehn, and Nicolai}}]{KKN09}
\bibinfo{author}{\bibfnamefont{A.}~\bibnamefont{Kleinschmidt}},
  \bibinfo{author}{\bibfnamefont{M.}~\bibnamefont{Koehn}}, \bibnamefont{and}
  \bibinfo{author}{\bibfnamefont{H.}~\bibnamefont{Nicolai}},
  \bibinfo{journal}{Phys Rev D} \textbf{\bibinfo{volume}{80}},
  \bibinfo{pages}{061701(R)} (\bibinfo{year}{2009}), \eprint{arXiv:0907.3048
  [gr-qc]}.

\bibitem[{\citenamefont{Koehn}(2011)}]{Koe11}
\bibinfo{author}{\bibfnamefont{M.}~\bibnamefont{Koehn}} (\bibinfo{year}{2011}),
  \eprint{arXiv:1107.6023v1 [gr-qc]}.

\bibitem[{\citenamefont{Humphreys}(1972)}]{Hum72}
\bibinfo{author}{\bibfnamefont{J.~E.} \bibnamefont{Humphreys}},
  \emph{\bibinfo{title}{Introduction to Lie Algebras and Representation
  Theory}} (\bibinfo{publisher}{Springer-Verlag New York},
  \bibinfo{year}{1972}).

\bibitem[{\citenamefont{Damour et~al.}(2001)\citenamefont{Damour, Henneaux,
  Julia, and Nicolai}}]{DHJN01}
\bibinfo{author}{\bibfnamefont{T.}~\bibnamefont{Damour}},
  \bibinfo{author}{\bibfnamefont{M.}~\bibnamefont{Henneaux}},
  \bibinfo{author}{\bibfnamefont{B.}~\bibnamefont{Julia}}, \bibnamefont{and}
  \bibinfo{author}{\bibfnamefont{H.}~\bibnamefont{Nicolai}},
  \bibinfo{journal}{Phys Lett} \textbf{\bibinfo{volume}{B509}},
  \bibinfo{pages}{323} (\bibinfo{year}{2001}), \eprint{arXiv:hep-th/0103094}.

\bibitem[{\citenamefont{Fuchs and Schweigert}(2003)}]{FS}
\bibinfo{author}{\bibfnamefont{J.}~\bibnamefont{Fuchs}} \bibnamefont{and}
  \bibinfo{author}{\bibfnamefont{C.}~\bibnamefont{Schweigert}},
  \emph{\bibinfo{title}{Symmetries, Lie Algebras and Representations}},
  Cambridge Monographs on Mathematical Physics (\bibinfo{publisher}{Cambridge
  University Press}, \bibinfo{year}{2003}).

\bibitem[{\citenamefont{Brion}(1988)}]{Bri88}
\bibinfo{author}{\bibfnamefont{M.}~\bibnamefont{Brion}}, \bibinfo{journal}{Ann
  Sci Ecole Norm Sup} \textbf{\bibinfo{volume}{21}}, \bibinfo{pages}{653}
  (\bibinfo{year}{1988}).

\bibitem[{\citenamefont{Baldoni et~al.}(2011)\citenamefont{Baldoni, Berline,
  De~Loera, Koeppe, and Vergne}}]{BBDLKV11}
\bibinfo{author}{\bibfnamefont{V.}~\bibnamefont{Baldoni}},
  \bibinfo{author}{\bibfnamefont{N.}~\bibnamefont{Berline}},
  \bibinfo{author}{\bibfnamefont{J.}~\bibnamefont{De~Loera}},
  \bibinfo{author}{\bibfnamefont{M.}~\bibnamefont{Koeppe}}, \bibnamefont{and}
  \bibinfo{author}{\bibfnamefont{M.}~\bibnamefont{Vergne}},
  \bibinfo{journal}{Math Comput} \textbf{\bibinfo{volume}{80}},
  \bibinfo{pages}{297} (\bibinfo{year}{2011}).

\end{thebibliography}

\end{document}